\documentclass{amsart}

\title{A Geometric Interpretation of  Ranicki Duality}
\author{Frank Connolly}

\address{Department of Mathematics, University of Notre Dame, Notre Dame IN 46556 U.S.A.}
\email{connolly.1@nd.edu }

\date{\today}

\usepackage{amsthm,amsmath,amssymb,eucal}
\usepackage{amscd,mathrsfs,verbatim,url}

\usepackage[colorlinks]{hyperref}

\newtheorem{thm}{Theorem}[section]

\newtheorem{cor}[thm]{Corollary}
\newtheorem{lem}[thm]{Lemma}

\newtheorem{nota}[thm]{Notation}
\newtheorem{prop}[thm]{Proposition}

\theoremstyle{definition}
\newtheorem{defn}[thm]{Definition}
\newtheorem{rem}[thm]{Remark}

\numberwithin{equation}{section}
\DeclareMathAlphabet{\matheurm}{U}{eur}{m}{n}

\newcommand{\Z}{\mathbb{Z}}

\newcommand{\s}{\sigma}

\newcommand{\cA}{\mathcal{A}}
\newcommand{\cB}{\mathcal{B}}

\newcommand{\op}{\mathrm{op}}

\newcommand{\p}{\partial}

\newcommand{\eps}{\varepsilon}

\newcommand{\lra}{\longrightarrow}

\newcommand{\rk}{{(R,K)}}

\begin{document}

\begin{abstract}
Consider a commutative ring $R$ and a simplicial map, $X\overset{\pi}{\lra}K,$ of finite simplicial complexes. The simplicial cochain complex of $X$ with $R$ coefficients, $\Delta^*X,$ then has the structure of an $(R,K)$ chain complex, in the sense of Ranicki \cite{ram}. Therefore it has a Ranicki-dual $ (R,K)$ chain complex,
$T \Delta^*X$. This (contravariant) duality functor  $T:\cB R_K\to \cB R_K$ was defined algebraically on the category of $(R,K)$ chain complexes and $(R,K)$ chain maps.

Our main theorem, \ref{MT}, provides a natural $(R,K)$ chain isomorphism:
\[
T\Delta^*X\cong C(X_K) 
\]
where $C(X_K) $ is the cellular chain complex of a (nonsimplicial) subdivision,  $X_K$,
of the complex $X$. The $(R,K)$ structure on $C(X_K) $
arises geometrically.
\end{abstract}

\dedicatory{To the memory of Andrew Ranicki}
\maketitle
\section{Introduction; Description of Results}
This article is an addition to a  theory of blocked surgery,  pioneered  by Ranicki and augmented by others in \cite{ram},\cite{rw},\cite{tmaf},\cite{dr},\cite{FJ},\cite{cdk},\cite{cdk2}. It is still in a developing state. It gives a new geometric interpretation of Ranicki's notion of the dual,
$TC$, of an $(R,K)$-chain complex $C$ when $C$ itself arises geometrically --
in particular when  $C=\Delta^*X$, if $(X,\pi)$ is a \emph{$K$-space} (defined in \ref{ksp}). This result, Theorem \ref{MT}, is the main theorem. It says, roughly, that
$T\Delta^*X = C(X_K)$, for the CW-complex $X_K$.

It is also our aim to give a transparent  definition of this duality functor $T$,  a clear treatment of
Ranicki's natural transformation $e:T^2\to id$. and  a  simple proof that $e_C:T^2C\to C$ is an $(R,K)$ chain equivalence for all  $C$.

Our   larger goal is to facilitate applications of Ranicki's theory to geometric questions such as the topological rigidity of non-positively curved groups
as in   \cite{FJ}, \cite{fj2}, \cite{cdk},\cite{cdk2} when those groups have elements of finite order.

The typical input of Ranicki's theory is a ``$K$-blocked normal map". By this we mean a  degree-one normal map of pl-manifolds, $(f, b): M^n\to X^n$, (as in \cite{br} or \cite{wa1}) plus  a ``control map", $\pi:X\to K$, where $K$ is a  finite simplicial complex. One seeks to understand the obstruction to obtaining a normal cobordism
 of $(f,b)$ to a ``$K$-blocked  homotopy equivalence", $N^n\to X$ (a map which is a homotopy equivalence over each ``block" in $K$).

 In the classical case ($K=$ point; \cite{br}, \cite{Wa2},\cite{wa1})  one has the ``surgery obstruction"  $\sigma(f,b)\in L_n(\Z[\pi_1(X)])$ to such a normal cobordism. 
 This functor $L_n()$, was generalized
 in \cite{ram} to yield obstruction groups $L_n(\cA)$ for any ``category-with-chain-involution" $ (A, *, \epsilon)$. Here $A$ is an additive category,  $\cB A\overset{*}{\to} \cB A$ is a contravariant functor satisfying certain conditions, on the category $\cB A$, of
 finite chain complexes in $A$, and $\epsilon:(*)^2\to id$, is an equivalence in the homotopy category of $\cB A$.
 
 %Elements of $L_n(\cA)$ are cobordism classes $[C,\psi]$ of ``$n$-dimensional quadratic chain complexes
% in $\cA$". Here $C\in Ob(B(A))$ is a  chain complex and $\psi$ (the ``quadratic structure on $C $ '') is  a chain-homotopy-analogue of a nonsingular quadratic form. 

 Ranicki, in \cite{ram},  starts with a finite complex $K$ and a category with chain involution, $\cA= (A, *, \eps)$ as above.
 He then constructs the additive category $A_K$ of $K$-blocked objects from $A$, and $K$-blocked
 $A$-maps. From $*$, and $\eps$, he defines the \emph{Ranicki Duality Functor}   $T:\cB(A_K)\to \cB(A_K)$,  and the  natural transformation   $e: T^2\to id_{\cB(A_K)}$. This construction allows one
  to define the surgery obstruction groups, $L_n(\cA_K)$ where $\cA_K =(A_K, T, e)$. 
 
 This applies directly to  a $K$-blocked normal map, $M^n\overset{(f,b)}{\lra} X^n\overset{\pi}{\to} K$. Here $A=A(R)$, the category of finitely generated free  modules over a fixed commutative ring $R$. 
 We  write $A R_K$ for $(A  R)_K$ and $\cB R_K$ for $\cB(A R_K)$. Its objects are \emph{$(R,K)$-chain complexes}. So  the
 simplicial cochain complexes of $X$ and $M$ denoted $\Delta^*X$ and $\Delta^*M$, and the simplicial chain complexes, $\Delta X'$ and $\Delta M'$, are $(R,K)$-chain complexes. (See Section \ref{sec3}). 
 
 However, Ranicki's original definition of $\cB R_K \overset{T}{\lra}\cB R_K$ is rather scattered and nonconceptual.
Indeed  his  assertion in \cite{ram} of the \emph{crucial} theorem that $(\cA R_K, T, e)$ is a category with chain involution
was only proved completely in  2018 (by Macko and Adams-Florou, \cite{tmaf}).

 \

This paper interprets Ranicki's notions geometrically. Section \ref{sec 2} fixes chain-complex conventions. Section \ref{sec3} reviews Ranicki's concepts concerning $(R,K)$ complexes while attempting to simplify notation.  In Section \ref{sec4} we introduce the $(R,K)$ chain complex $C\otimes_KD$, defined if $D$ is an $(R,K)$ complex and $C$ is an $(R,K^{op})$ complex.  This  complex $C\otimes_KD$ is a certain quotient of $C\otimes_RD$.

Our definition (see \ref{rad}) of the Ranicki dual $TC$, of an $(R,K)$ complex $C$, is:
\[
TC=C^*\otimes_K \Delta^*K.
\]

  In Section \ref{bc} we show, using work of M. Cohen \cite{cohen},  that each $K$-space $(X,\pi)$ defines a certain regular CW-complex  $X_K$, whose cellular chain complex has a natural $(R, K)$ structure.  Therefore
from each $K$-space $(X,\pi)$ we obtain three $(R,K)$ chain complexes:
\begin{enumerate}
\item $\Delta^*X$, the simplicial cochain complex of $X$ (Definition \ref{ksp}).
\item $C(X_K)$, the cellular chain complex of the CW complex $X_K$ (Section \ref{isom}).
\item $\Delta X'$, the simplicial chain complex of $X'$, the barycentric subdivision of $X$ (Definition \ref{derived}).
\end{enumerate}
This paper shows that these three are closely related by $T$. Our
main result, Theorem \ref{MT}, exhibits an isomorphism of $(R,K)$ chain complexes: 
\[
\Phi_X: T\Delta^*X\cong C(X_K)
\]
When $X$ is a pl-manifold, Poincare duality then becomes an $n$-cycle in the $(R,K^{\op})$ complex, $Hom_{(R,K)}(TC(X_K), C(X_K))$.

We then use the work of McCrory, \cite{McC}, to prove there are $(R,K)$ chain homotopy equivalences:
$\qquad
T\Delta X'\simeq \Delta^*X;\quad \quad C(X_K)\simeq \Delta X'.
$

  This regular CW complex $X_K$ is a subdivision of $X$. The derived complex $X'$,  is  a simplicial subdivision of $X_K$. 
  In fact, for each simplex $T$ of $X$ and each face $\s$ of $\pi(T)\in K$,   there is a single cell $T_\s$ of $X_K$. This cell $T_\s$ turns out to be simply $(\pi\mid T)^{-1}|D(\s,\pi(T))|,$ where $D(\s,\pi(T))$ is the dual cell of $\s$ in $\pi(T)$.

\

The author is indebted to Jim Davis for his helpful comments on this paper.
  
 % In  geometrically defined cases ( when $C=\Delta^*X$, for any $K$-space $X\overset{\pi}{\lra} K$)  we then give a \emph{geometric interpretation}
  %of $T\Delta^*X$.  We define a regular CW complex,  $D_KX$ whose cellular chain complex, $C(D_KX)$ is naturally an $(R,K)$-chain complex.  We then exhibit a natural $(R,K)$ chain isomorphism
 % \[
 % \Phi_X: T\Delta^*X\;\cong \;C(D_K X).
  %\]

 \section{ Chain Complex Conventions}\label{sec 2} 
Throughout  this paper, $R$ denotes a fixed commutative ring; $AR$ is the additive category of finitely generated free $R$ modules.
     
For any additve category $A$ we will write $\cB A$ for the additive category of finite chain complexes, $C=\{C_q, \p_q\}_{q\in \Z}$ and chain maps $f=\{f_q:C_q\to D_q\}_{q\in\Z}$ from $A$.
(\emph{Finite} means: $C_q=0$ for all but finitely many $q$). We abbreviate $\cB (AR)$ to $\cB R$.

As usual  two chain maps $f,g: C\to D $ are chain homotopic if there  is a sequence of $A$ maps,  $h=\{h_q:C_q\to D_{q+1}\}$, for which
$d^D_{q+1}h_q + h_{q-1}d^C_q = g_q-f_q \quad \forall q$.

We regard $A$ as the full subcategory of $\cB A$ consisting of chain complexes concentrated in degree zero.

Let $C,D\in Ob(\cB R)$. The complexes $C\otimes_R D$,  and $Hom_R(C,D)$ in $Ob(\cB R)$, are:
\begin{gather}
 \boxed{(C\otimes_R D)_q = \sum_{r\in\Z} C_r\otimes_R D_{q-r};\qquad  Hom_R(C,D)_q =\sum_{r\in \Z}Hom_R(C_r, D_{q+r}) } \;\text{ and:}\notag \\
   \boxed{ 
 d^{C\otimes D}(x\otimes y)= d^Cx\otimes y +(-1)^{|x|} x\otimes d^Dy;\qquad  d^{Hom}\phi =d^D\circ \phi - (-1)^{|\phi| } \phi\circ d^C }\notag
\end{gather}
The \emph{evaluation map}, $ eval_{C,D}:Hom_R(C,D)\otimes_R C\lra D$ is the \emph{R-chain map: }
\[
 eval_{C,D}(f\otimes x)= f(x).
 \]

Write $ev_C:C^*\otimes C\to R$ for $eval_{C,R}$.

 The contravariant functor  $\cB R\overset{*}{\to} \cB R$ is : $C^*= Hom(C,R);  \; f^* = Hom(f, 1_R)$.
 
 Therefore we have:
 \[
 (C^*)_{-q} = Hom_R(C_q,R); \qquad  d^{C^*}_{-q} = (-1)^{q+1}(d^C_{q+1})^*:(C^*)_{-q}\to (C^*)_{-q-1}.
 \] 
 
The functor $*$ comes with a natural equivalence, $\eps: (*)^2\to 1_{\cB R}$. Specifically,  the chain isomorphism
$\epsilon_C: C^{**}\to C$ is  characterized by the identity:
\[
 a(\eps_C(\alpha)) = (-1)^q \alpha(a) \qquad \forall \alpha\in (C^{**})_{q},\; a\in (C^*)_{-q} .
\]

 \section{ Basic Definitions For $(R,K)$ Chain Complexes}\label{sec3}
 \begin{defn} Let $K$ be a finite poset with partial order  $\leq$.  For example, K can be a finite  simplicial complex where $\s\leq\tau$ means $\s$ is a face of $\tau$.
  $K^{op}$ denotes the same set with the opposite partial order. 
  \begin{enumerate}
  \item An \emph{(R,K) module} is an ordered pair 
 $M= (M(K), \{M(\sigma)\}_{\sigma \in K} \,)$ such that: 
 \begin{enumerate} 
\item  $ M(K)$ and each $M(\s)$ are $R$-modules in $Ob(A R)$;
\item $M(K) = \oplus_{\sigma\in K} M(\sigma)$.
\end{enumerate}
More generally, for any $S\subset K$ we write: $M(S)= \oplus_{\sigma\in S} M(\sigma).$  
    \item An \emph{$(R,K)$ map  $M\overset{f}{\to}N$} of $(R,K)$ modules is a map $M(K)\overset{f}{\to} N(K)$ of \quad $R$ modules, whose components,  $f(\tau,\sigma):M(\sigma)\to N(\tau),$ satisfy:
 \[ f(\tau,\sigma)=\,0  \text{ unless } \tau\geq\s.
 \]
 \item The additive category of $(R,K)$ maps and modules is written $AR_K$. 
 
 We abbreviate the category of chain complexes, $\cB(AR_K)$,  to $\cB R_K$.
 
 \item An object $C =\{ C_q, \p_q \}_{q\in \Z}$ of  $\cB R_K$ is an \emph{($R,K)$ chain complex}. We then write $C(K)$ for $\{C_p(K), \p_p\}_{p\in \Z}$, an $R$-chain complex in $ob(\cB R )$.  
 
 Note:  $C\in ob(\cB R_K)$
 is  specified by specifying the $R$ complex $C(K) $ and the required collection $\{C_q(\s)\}_{\s\in K, q\in\Z} $\; of $R$ submodules.

  \item 
  Let $C,D\in ob(\cB R_K).$ $Hom_{(R,K)}(C,D)$ is the $(R,K^{op})$ complex such that:
  \begin{enumerate}
  \item $Hom_{(R,K)}(C,D)(K)$ is the  subcomplex  of $Hom_R(C(K), D(K)) $ given by those $f =\{f_q : C_q\to D_{q+|f|}\}_{q \in \Z }$ for which each $f_q $ is an $(R,K)$ map.
  \item $Hom_{(R,K)}(C, D)_p(\s)$ is the set of  $f\in Hom_{(R,K)}(C, D)(K)_p$ satisfying: 
  \[
  f_q \mid_{C_q(\tau) } =0 \text{ if } \tau \neq \s,\; \forall\; q.
  \]
  \end{enumerate}
 \item We say a sequence of chain maps  $0\to C'\overset{i}{\to} C\overset{j}{\to} C''\to 0$ in $\cB R_K$  is \emph{exact} if for each $\s\neq \tau$, $i(\s, \tau)=0,\; j(\s,\tau)=0, $ and, for all $q$,  the corresponding sequence,  $0\to C'_q(\s)\to C_q(\s)\to C''_q(\s)\to 0$. is an exact sequence in $\cA R$. 
We then say $i$ is an $(R,K)$ \emph{monomorphism} and $j$ is an $(R,K)$ \emph{epimorphism}.

\item Note that $*$ specifies a contravariant functor, $\cB R_K\overset{*}{\lra}\cB R_{K^{op}}$,  provided that we define\;
$(C^*)_q (\s) $ as $(C_{-q}(\s))^* $ and $d^{C^*} $  as $d^{C(K)^*} $ for $C\in ob(\cB R_K)$. 
$\cB R_K\overset{*}{\lra}\cB R_{K^{op}}$
preserves exactness and  homotopy. The transformation $\eps_C: C^{**}\to C$ of Section \ref{sec 2} is an $(R,K)$ isomorphism, for  all $C\in ob(\cB R_K).$ 

\item We say $S\subset K$ is \emph{full}   if, whenever $\rho,   \tau \in K$, then
 $\{\s\mid \rho\leq \s\leq \tau\} \subset S$.  \newline
 Let $C$ be an $(R.K)$ complex. Let $S$ be full.  Define  $\p_q^{C(S)}:C_q(S) \to C_{q-1}(S)$ by:    $\p_q^{C(S)} x = \sum_{\tau\in S}\p^C(\tau, \s)x$ for all $\tau, \s \in S$ and all $x\in C(\s)$.
 
  Then $C(S):=\{ C_q(S), \p_q^{C(S)}\}_{q \in \Z}$ is an $R$ chain complex.

\item In particular, a simplicial complex $K$ is a poset. As usual $\Delta_*(K;R)=\{ \Delta_q(K;R),\p_q\}_{q\in \Z}$ denotes the simplicial chain 
complex of $K$. One can choose a basis, $bK$ for $\Delta_*(K;R)$ consisting of one oriented $q$-simplex, $\s=$\newline $\langle v_0\dots, v_q\rangle\in \Delta_q(K;R)$ for each $q$-simplex with vertices $v_0,\dots v_q$, of $K$. Recall:
$\langle v_0,\dots,v_q\rangle = sgn(\pi)\langle v_{\pi(0)},\dots,v_{\pi(q)}\rangle$
for each $\pi\in S_{p+1}$. 
The oriented $q$-simplex $\s\in\Delta_q(K;R)$ defines a dual cochain $\s^*\in \Delta^*(K;R)_{-q}$ such that $\s^*(\tau) =0$ for all $ \tau\neq \pm \s$, and $\s^*(\s)=1$. 

One then defines $\s^{**}\in \Delta_q(K;R)^{**}$
 by: $\eps(\s^{**}) = \s$. 
 
 Here $\Delta^*(K;R)$ is the simplicial cochain complex of $K$.
 
 Each simplex $\s\in K$ defines  subcomplexes, $\overline{\s}$ and $ \p\s$, 
and a subset $st(\s)$:
\[
\qquad\quad \overline{\s}=\{\tau\in K\mid \tau\leq \s\}; \quad \p\s=\{\tau\in K\mid \tau< \s\}; \quad st(\s)=\{\tau\in K\mid \tau\geq \s\}
\] 
\item The incidence number $[\tau, \s]\in \{1,-1, 0\}$ is defined for any oriented
simplices $\s, \tau$ of $K$. It satisfies: $\p_q(\s)= \sum_{\tau\in bK} [\s, \tau] \tau$ for any basis, $bK$ of oriented simplices of $K$.  $[\s,\tau]\neq 0$ iff $\tau$ is a codimension-one face of $\s$.

 \end{enumerate}
\end{defn}
\

\begin{defn}\label{ksp} ($K$-spaces, $\Delta^*X$ and $\Delta X$) 

Let $K$ be a finite simplicial complex. 
A \emph{$K$-space} is a pair $(X, \pi)$ where $X$ is a finite simplicial complex and $|X|\overset{\pi}{\to} |K|$ is a simplicial map, $X\to K$. 
A \emph{map of
$K$-spaces}, $(X,\pi_X)\to (Y, \pi_Y)$ is a simplicial map $f:|X|\to |Y|$ satisfying: $\pi_Y f=\pi_X$.

 Let $(X,\pi)$ be a $K$-space.
  
 $\Delta X$ denotes the $(R,K^{op})$ complex for which $\Delta X(K)= \Delta_*(X;R)$.
 For each $\s\in K$, $(\Delta X)_p(\s)$ is the submodule generated by  oriented $p$-simplices in $ \Delta_p(X;R)$ whose underlying p-simplex, $S\in X$,  satisfies $\s= \pi (S)\in K$.   
  
 By definition, $\Delta^*X= *(\Delta X)$. Therefore $\Delta^*X(K) = Hom_R(\Delta_*(X;R), R)=\Delta^*(X; R)$, the simplicial cochain complex of $X$. For each $\s\in K$, $(\Delta^*X)_{-p}(\s)$ is therefore the submodule spanned by all $S^*$ for which $S\in\Delta_p(X;R)$ is an oriented simplex  and $\s=\pi(S)\in K$.

 A map $f:X \to Y$ of $K$-spaces induces an $(R,K)$ chain map $f^*:\Delta^*Y\to \Delta^*X$ and an $(R,K^{op})$ chain map $f_*:\Delta Y\to \Delta  X$. 
\end{defn}

The next lemma will be used in section \ref{RD}.
\begin{lem} \label{clem}Suppose $S\in K$ and there is no \,$ \tau\in K$ for which $S<\tau$. The $K$-space $(\overline{S},inclusion)$ specifies the 
$(R,K)$ complex $\Delta^*\overline{S}$.  
Then $\Delta^*\overline{S}(st(\s)) $ is  a contractible $R$-complex for all $\s\in K$ such that $\s\neq S$. Also $\Delta^*\overline{S}(st(S)) = R S^*.$
\end{lem}
\begin{proof}It is obvious that $\Delta^*\overline{S}(st(S)) = R S^*$ (after orienting $S$) and that
$\Delta^*\overline{S}(st(\s)) =0$ if $\s$ is not a face of $S$. So we assume
$\s<S$. Let $\tau$ be the complementary face of $\s$ in $S$. Then the joins,
$\overline{S} = \s* \tau$ and $\p\s *\tau$ are contractible simplicial complexes. Note $st(\s)= \overline{S} -\p\s *\tau$. Consequently, $\Delta^*\overline{S}(st(\s)) = \Delta^*(\s * \tau, \p\s * \tau;R)$ is a contractible chain complex.
\end{proof}

\

\section{ $C\otimes_K D$ and the isomorphism $ Hom_{(R,K)}(D, C^*) \cong  (C\otimes_K D)^*$  }\label{sec4}
\emph{ For the rest of this paper, $K$ denotes a finite simplicial complex. }

\emph{Moreover, throughout this section,
 $C$ denotes an
$(R,K^{op})$ complex and   $D$ denotes an $(R,K)$ complex.}

 \

In $K$, the star of any simplex ,  $st(\s)$, as well as $K-st(\s)$ are full in $K$. Moreover the chain complex $C(K-st(\s))$ is a subcomplex of $C(K)$ and the resulting inclusion fits into a short exact sequence of chain maps in $\cB R$:
\[
    0\to C(K-st(\s))\overset{i_\s}{\lra} C(K)\overset{p_{st(\s)}}{\lra}C(st(\s)) \to 0
\]
Here $C(K)\overset{p_{st(\s)}}{\lra}C(st(\s))$ is defined by: $p_{st(\s)}|_{C_q(st(\s))}= 1_ {C_q(st(\s))}$;\;  and \newline $p_{st(\s)}|_{C(K-st(\s))}=0$.

\begin{defn}\label{ckd} ( $C\otimes_KD$,\quad   $C\otimes_RD$,\quad and $C\otimes_RD \overset{\pi_{C,D}}{\lra}C\otimes_K D$). 

Let $C$ be an
$(R,K^{op})$ complex and $D$ be an $(R,K)$ complex. 
 \begin{enumerate}
\item Let $C\otimes_RD$ be the $(R,K)$ complex for which:
\[
(C\otimes_RD)(K)= C(K)\otimes_R D(K); \quad (C\otimes_R D)_q(\rho)= (C(K)\otimes_R D(\rho))_q \quad \forall \rho\in K, q\in \Z
\]
\item

Let  $C\otimes_KD$ be the $(R,K)$ complex for which:
\begin{enumerate}
\item \quad $(C\otimes_K D)_q(K)= \sum_{\rho\in K} (C(st(\rho))\otimes_R D(\rho))_q       \quad \forall q\in\Z      $   
\item  \quad $(C\otimes_K D)(\rho)= C(st(\rho))\otimes_R D(\rho)      \quad \forall \;\rho\in K $
\item The map $ C\otimes_RD\overset{\pi_{C,D}}{\lra} C\otimes_KD $ is an $(R,K)$ chain epimorphism, if we define $\pi_{C,D}$ by requiring that 
$\pi_{C,D}(\s,\rho)=0$ for $\s\neq \rho$  and
\[
\pi_{C,D}(\rho, \rho)= p_{st(\rho)}\otimes_R 1_{D(\rho) }:C(K)\otimes_R D(\rho) \lra C(st(\tau))\otimes _R D(\rho).
\]
 \end{enumerate}
 \end{enumerate}
Expicitly,  for any $\rho\leq \tau$ and $x\otimes_Ry\in C_r(\tau)\otimes_RD_{q-r}(\rho)\subset 
(C\otimes_KD)_q(\rho)$, we have

\begin{equation}
\boxed{\label{d_K}
d^{C\otimes_K D}(x\otimes y) = \sum_{\{\s| \rho\leq \s\leq\tau\} }
  d^C(\s,\tau)x \otimes y + (-1)^r      x\otimes  d^D(\s, \rho) y.
}
\end{equation}
\end{defn}

\

We now show that $(C\otimes_KD)^*$ is a convenient expression for $Hom_{(R,K)}(D, C^*)$:
 \begin{lem} \label{nate}There is a natural isomorphism $\Psi$ of  functors, denoted,
 \begin{equation}
\Psi_{C, D}: Hom_\rk (D, C^*)\cong (C\otimes_K D)^*
\end{equation}
for any $(C,D)\in Ob(\cB R_{K^{op}}\times \cB R_{K})$. 
 \end{lem}

\begin{proof} Suppose $f$ is in $Hom_{(R,K)}(D, C^*)_q(\s)$ for some $\s\in K$ and $q\in\Z$. Define an $R$-map, $\Psi(f): C(st(\s))\otimes D(\s)_{-q}\to R$, by the formula:
\[
\Psi(f)(x\otimes y)= (-1)^{|x| |y|} f(y)(x) \qquad \text{for } x\otimes y\in (C\otimes_KD)_{-q}(\s). 
\]
The same formula yields $0$,  if $x\otimes y$ is in  $(C\otimes_KD)(\tau)_{-q}$ for $\tau\neq \s$. One easily sees that this rule (i.e. $ f\mapsto \Psi(f) $   gives an
isomorphism,
\[
\Psi_{C,D}: Hom_{R K}(D,C^*)\overset{\cong}{\lra} (C\otimes_K D)^*
\]
of $(R,K^{op})$ complexes for all such $C, D$. Naturality is straightforward.
\end{proof}

%\begin{defn}

%\end{defn} 

\section{Ranicki Duality and the $(R,K)$ chain equivalence $e:T^2\to 1_{\cB R_K} $ }\label{RD}

\begin{defn} \label{rad} Ranicki Duality is the contravariant functor $\cB R_K\overset{T}{\lra}\cB R_K$ \newline
defined for a chain complex $C\in Ob(\cB R_K)$ and a $(R,K)$ chain map, $f:C\to D$ by:
\[
TC= C^*\otimes_K\Delta^*K\qquad \qquad Tf= f^*\otimes_K  1_{\Delta^*K}
\]
\end{defn}
$\Delta^*K$ comes from the $K$-space, $(K,1_K)$.
After examining \cite{ram}, p. 75 and p.26, lines -6 to -4 one can see that
this is in agreement with the definition  indicated there, up to isomorphism and differences in sign conventions. In particular compare our formula for $d^{C\otimes_KD}$ with that on p.26, line -5 of \cite{ram}.
\begin{cor} $T$ is an exact homotopy functor.
\end{cor}
\begin{proof}: By Lemma \ref{nate}, $TC= C^*\otimes_K \Delta^* K$ is isomorphic to $Hom_{(R,K)}(\Delta^*K, C)^*$ (since $\eps_C:C^{**}\cong C$ for all $C$).
But $C\mapsto C^*$ and $C\mapsto Hom(\Delta^*K, C)$ are both exact homotpy functors. The result follows.   \end{proof}
 
 We now want to show that $T^2C$ and $C$ are $(R,K)$-chain equivalent. See \ref{cheq}.

\begin{defn}(of  $E_C:Hom_{(R,K)}(\Delta^*K, C)\otimes_K 
\Delta^*K\to C$).

Let $C$ be an $(R,K)$ complex. 

Consider the evaluation chain map, $Eval_{A,B}:Hom_R(A,B)\otimes_R A\to B$, when $A=\Delta^*K(K)$ and
$B=C(K)$. Its restriction to $(Hom_{(R,K)}(\Delta^*K, C)\otimes_R\Delta^*K)(K),$ denoted $E'_C,$
is an $(R,K)$ chain map, 
\[
E'_C:Hom_{(R,K)}(\Delta^*K, C)\otimes_R\Delta^*K\to C
\]
 (by definition of an $(R,K)$ map). Moreover, for each $\s\in K$,  $E'_C$ annihilates $Hom_{(R,K)}(\Delta^*K, C)(K-st(\s))\otimes_R\Delta^*K(\s)$. Therefore $E'_C$ descends uniquely  to an $(R,K)$ chain map, 
\[
E_C:Hom_{(R,K)}(\Delta^*K, C)\otimes_K\Delta^*K\to C,\quad\qquad E_C(f\otimes \s^*)=f(\s^*).
\]
satisfying:  $E'_C= E_C\circ\pi_{H, \,\Delta^*K} $ . 
Here $H= Hom_{(R,K)}(\Delta^*K, C)$
(see \ref{ckd}).

$E$ is obviously natural in $C$.
\end{defn}
For each $(R,K^{op})$ complex $C$, set $\Psi_C= \Psi_{C, \Delta^*K}$

In view of Lemma \ref{nate}.  we have an $(R, K)$ chain isomorphism:
\[
\Psi_C\otimes 1_{\Delta^*K}: Hom_{(R,K)}(\Delta^*K, C)\otimes_K\Delta^*K\overset{\cong}{\lra}(C^*\otimes_K\Delta^*K)^*\otimes_K \Delta^*K= T^2C .
\]

 \begin{defn}For each $(R,K)$ complex $C$ define $e_C:T^2C\to C$ as the unique map for which $E_C= e_C\circ (\Psi_C\otimes 1_{\Delta^*K})$.
 
 Note $e_C$ is an $(R,K)$ chain epimorphism and $e$ is a natural transformation.
 \end{defn}

\begin{thm} \label{cheq} $e_C: T^2C\overset{e_C}{\lra} C$ is an $(R,K)$ chain equivalence, for each $(R,K)$ complex $ C$.

\end{thm}
\begin{proof} By \cite{ram} (Proposition 4.7), we  need only prove that $e_C(\s,\s): T^2C(\s) \to C(\s)$ is an $R$-chain equivalence, for all $\s\in K$. 

Case I: Assume there is a simplex $S\in K$ for which: $C(\s)=0 \quad \forall \;\s\neq S$. 

We need only show $e_C(S,S)$ is a chain isomorphism, and
$T^2C(\s)$ is contractible for $\s\neq S$.
We compute, for all $\s\in K$, in view of the restriction on $C$:
\begin{align*}
     TC(st(\s)) &= (C^*\otimes_K \Delta^*K)(st(\s)) \\
     =(C^*\otimes_R \Delta^*\overline{S})(st(\s)) &=C^*(S)\otimes_R\Delta^*\overline{S} (st(\s))\\    
   \text{So: }  T^2C(\s) \cong & C^{**}(S)\otimes_R \Delta^{**}\overline{S} (st(\s)) \otimes_R R\s^*
     \end{align*}
So for $\s \neq S, \; T^2C(\s)$ is contractible because  $\Delta^{**}\overline{S} (st(\s)) $ is contractible by \ref{clem}. In addition, after orienting $S$ we conclude
       \begin{align*}
     T^2C(S) = & C^{**}(S)\otimes_R R S^{**} \otimes_R R S^*,\\
     e_C(S,S)( c&\otimes S^{**}\otimes S^*)=  \pm \eps_C(c)\quad \forall \;c \text{ in } C^{**}(S).
     \end{align*}
     So $e_C(\s,\s)$ is a chain isomorphism for $\s= S$ and  a   chain equivalence for $\s\neq S$.
     This completes the proof in Case I.
     
     Case II (the general case): For any $C\neq 0$  in $\cB R_K$ one can choose some $S\in K$ for which $C(S)\neq 0$, and an exact sequence $0\to C'\overset{i}{\lra} C\overset{j}{\lra} C''\to 0$ for which $i(S,S):C'(S)\to C(S)$ is an
     isomorphism, and $C'(\s)= 0$ for $\s \neq S$. For example, choose $S$ to be of maximum dimension among $\{\s\in K \mid\; C(\s)\neq 0\}$).  
     
     The argument is by induction on the number 
      $n$, of $\s\in K$, for which $C(\s)\neq 0$.
     
      If $n=1$, Case I applies.  
     If $n>1$, by induction, $e_{C''}(\s,\s)$ and $e_{C'}(\s,\s)$ are \quad\quad
     $R$ chain equivalences. Also the commuting diagram below has exact rows.
  \[
   \begin{CD}
   0 @>>>C'(\s) @>>> C(\s)  @>>> C''(\s) @>>> 0\\
   @.  @Ve_{C'}(\s,\s)VV @Ve_{C}(\s,\s)VV @VVe_{C''}(\s,\s)V\\
   0@>>> T^2C'(\s) @>>> T^2C(\s) @>>> T^2C'' (\s) @>>> 0\\
   \end{CD}  
     \]
      Therefore $ e_C(\s,\s)$ is an $R$-chain equivalence for all $\s$. This completes the proof.
      \end{proof}
Note: The first proof of the above theorem appeared in \cite{tmaf}.

\

\section{ Construction of the Ball Complex $X_K$}\label{bc}

The purpose of this section is to construct the complex $X_K$ advertised in the introduction and establish
its properties.

\begin{nota} 
Let $X$  be a finite simplicial complex in a euclidean space, with vertex set $V_X$. Its underlying polyhedron is:  $|X|=\cup\{\sigma\;| \;\sigma\in 
X\}$. For each $p\geq 0$, $X_p$ denotes the set of $p$-simplices of $X$. 

 If $|X|$ is pl-homeomorphic to $I^n$  we say $|X|$ is a pl $n$-ball and write $\p X$ for the subcomplex for which
$|\p X|= \p |X|$.

Each p-simplex $\sigma \in X$ is the convex hull,  $[v_0, v_1,\dots, v_p]$,  of its  vertices in $V_X$. Its \emph{barycenter}  is $\hat{\sigma}:= \frac{1}{p+1} \sum_{i=0}^p v_i \in\sigma^\circ$. 

%The barycentric coordinate function of $v\in V_K$ is written
%$a_v: |K|\to [0,1]$.

Choose a point $b\sigma\in \sigma^\circ $ for each $\sigma\in X$. The derived complex $X'$ is defined as the unique simplicial subdivision of $X$ for which $V_{X'} = \{b\sigma \mid\; \sigma\in X\}$. 
$X'$ has one \newline $p$-simplex, $ [b\sigma_0, b\sigma_{1} \dots b\sigma_p]$,   for each decreasing sequence  $\sigma_0>\dots>\sigma_p$ of simplices of $X$. The ordered $p+1$ tuple
$ (b\sigma_0, b\sigma_{1}, \dots, b\sigma_p)$  specifies an oriented $p$-simplex denoted $\langle b\sigma_0, b\sigma_{1} \dots b\sigma_p\rangle\in \Delta_p(X';R)$ which we abbreviate to $\langle \sigma_0, \sigma_{1} \dots, \sigma_p\rangle$. These form a canonical basis for $\Delta_p(X';R)$ (in contrast to $\Delta_p(X;R)$).
\end{nota}
Because we want to use the McCrory cap product, we follow the orderings of \cite{McC} regarding simplices of $X'$.

\begin{defn} \label{derived} Let $(X,\pi)$
be a  \emph{$K$-space}. 
The \emph{derived complexes of $(X,\pi)$} are the simplicial subdivisions $X'$ of $X$,  and $K'$ of $K$ whose vertex sets $\{b\s \mid \s\in K\}$ and $\{bS\mid S\in X\}$ are chosen as follows:
\[
\text{If }   \sigma\in K, \qquad b\sigma:= \hat{\sigma}\in\sigma^\circ;
\]
\[
\text{If }  S\in X\text{ and } \sigma=\pi(S), \qquad bS:=  \text{ centroid of } (S\cap \pi^{-1} (\hat{\sigma}))\in S^\circ.
\]
By construction, $\pi(V_{X'})\subset V_{K'}$.  So $\pi$ is also a simplicial map from $X'$ to  $K'$, because $\pi$ is linear on each simplex of $X'$.

\

$X'$ provides a second geometric example, $\Delta X'$, of an 
$(R,K)$  complex:

 We define $\Delta X'$ by,
  \begin{enumerate}
  \item $\Delta X'(K)= \Delta_*(X';R)$.
  \item For each $\s\in K, p\in \Z,\;\; (\Delta X')_p(\s)$ is the submodule of $\Delta_p(X';R)$ spanned by all $\langle Q^0,\dots Q^p\rangle $ in $X'$ for which $\s= \pi(Q^p)$.
  \end{enumerate}
  It is straightforward to see that $\Delta X'$ is an $(R,K)$ complex.
\end{defn}
 
  \
 
The dual cone of a simplex $\sigma\in K$, denoted $D(\sigma, K)$, was first defined in \cite{poi}, Section 7. It is a subcomplex of $K'$ (and a pl ball if $K$ is a pl-manifold).  It gives rise to several ``dual" subcomplexes in $K' $ and $X'$ which we define now.
%Many of the combinatorics of this and other closely related dual subcomplexes were worked out by
%M.Cohen \cite{cohen}, whose results we use here. 
 
\begin{defn} Let $(X,\pi)$ be a $K$-space. Suppose $\s, \tau \in K$, $T\in X$.
\begin{enumerate}
\item $D(\sigma, K):=\{\langle \sigma_0, \sigma_1,\dots, \sigma_p\rangle\in K' \mid  \sigma\leq \sigma_p \}$
%\item $\dot{D}(\sigma, K):=\{\langle \sigma_0, \sigma_1,\dots, \sigma_p\rangle\in K' \mid  \sigma < \sigma_0\}$
\item $D(\sigma, \tau):\;=\,\{\langle \sigma_0, \sigma_1,\dots, \sigma_p\rangle\in K' \mid  \sigma\leq \sigma_p, \;\sigma_0 \leq \tau\}$
%\item $D_\sigma \pi:= \{\langle S_0, S_1,\dots, S_p\rangle \in X' \mid \sigma\leq \pi(S_0\}$
%\item  $\dot{D}_\sigma \pi:= \{\langle S_0, S_1,\dots, S_p\rangle \in X'\, \mid \sigma<\pi(S_0\}$
\item $D_\sigma T:=\{\langle S_0, S_1,\dots, S_p\rangle \in X' \mid \sigma\leq \pi(S_p),\; S_0\leq T\}$
\item $T_\sigma:= |D_\sigma T|$.
\end{enumerate}
\end{defn}
Of course,  $D(\sigma, \tau)=\emptyset$ unless $\sigma\leq \tau$, and  $ D_\sigma T=\emptyset$ unless $\sigma\leq \pi(T)$.

 %Note that
 % $|D_\sigma \pi|\;= \pi^{-1} (|D(\sigma, K)|)$, and 
  %$|\dot{D}_\sigma \pi|= \pi^{-1} (|\dot{D}(\sigma, K)|)$.

 \;$D_\sigma T$ is a subcomplex of $X'$. $D(\s , K) $   and $D(\s, \tau)$ are subcomplexes of $K'$.

\

%Let $B$ be a simplicial complex for which $|B|$ is pl homeomorphic to $I^n$ for some $n\geq 0$. We say the polyhedron $|B|$ is a pl n-ball. We write $\p B$ for the subcomplex of $B$ for which $\p |B| = |\p B|$.

\begin{lem}\label{ballem}Let $(X,\pi)$ be a $K$-space. Suppose $\sigma\in K$, $T\in X$, and $\sigma\leq \pi(T)$.  \begin{enumerate} 
\item $T_\sigma=|D_\sigma T|$ is a pl ball.  \; $dim(T_\sigma) = dim(T)-dim(\sigma)$.  

\item  $\p D_\sigma T=\p^i D_\sigma T \cup \p^o D_\sigma T$, 
(the inner and outer boundaries) where:
\[
\p^i D_\sigma T=\cup \{D_\rho T\;\mid \sigma< \rho \};\qquad \p^oD_\sigma T=\cup \{D_\sigma S\;\mid  S< T\}
\]
\item Suppose $\sigma<\pi(T)$.  Then $|\p^iD_\sigma T|$ and $ |\p^o D_\sigma T|$ are  pl balls of dimension
$dim(D_\sigma T)-1$, and  
\[\p (\p^iD_\sigma T)=\p(\p^oD_\sigma T)=\p^iD_\sigma T\cap \p^oD_\sigma T\;.
\]
\end{enumerate}
\end{lem}
\begin{proof} of (1):
For each vertex $v$ of $\tau$ note that, 
\[
|D(v,\tau)|=\{x\in \tau \;| \;a_v(x)\geq a_w(x),\text{ for all vertices $w$ of } \tau\}.
\]
where $a_v: |K|\to [0,1]$ denotes the barycentric coordinate function defined by  the vertex $v$.
This is a convex subset of   $\tau$.  So 
\[
 \; |D(\sigma, \tau)|=\cap_{v\in V(K)}  |D(v, \tau)|
\]
is also convex. Therefore $T_\s=(\pi_{|T})^{-1} (|D(\sigma, \tau)|)$ is also convex since  $\pi_{|T}: T\to \tau$ is  simplicial. So $T_\sigma  $ is a compact convex
polyhedron and therefore a pl ball.

 Since $|D(\sigma, \tau)|\cap \tau^\circ\neq \emptyset$, this operator $(\pi_{|T})^{-1} $ preserves codimension:
\[
dim(\tau) - dim(D(\sigma, \tau)) = dim (T) - dim (D_\sigma T).
\]
Since $dim(D(\sigma, \tau)) = dim(\tau)- dim(\sigma)$, we get: $dim(D_\sigma T)= dim(T)-dim(\sigma)$.
\end{proof}
\begin{proof} of (2): See \cite{cohen}, Proposition 5.6(2), applied to   $\pi_{|\overline{T}} :\overline{T}\to \pi(\overline{T})$.
\end{proof}
\begin{proof} of (3): The equation in (3), and the fact that  $|\p^iD_\sigma T|$ and $ |\p^o D_\sigma T|$ are both pl manifolds, are proved in \cite{cohen} (Proposition 5.6 (3),(4)). To show $|\p^iD_\sigma T|$ is a pl ball, it suffices to note that
it collapses to the vertex $bT$, and so  $|\p^iD_\sigma T|$ is a regular neighborhood of $bT$ in $|\p D_\sigma T|$ (by 3.30 of \cite{rs}).  Then by  3.13 of\cite{rs},  $\p^oD_\sigma T$ is also a pl ball.
\end{proof}

\begin{defn} (\cite{rs} p.27) A \emph{ball complex} is a finite   collection $ Z=\{B_i\}_{i\in I}$ of pl balls in  
a euclidean space, such that each point of $|Z|:=\cup \{ B\;\mid\; B\in Z\}$ lies in the interior of precisely one ball of $Z$, and the boundary of each $B\in Z$ is a union of balls  of lesser dimension of $Z$. Therefore $(|Z|, Z)$ is a regular CW-complex. 

Let $Z$ and $Y$ be ball complexes A pl map $f:|Z|\to |Y|$ is a \emph{ map of ball complexes } if for each ball $B$ of $Z$,  $f(B)$ is a ball of $Y$.
\end{defn}

\begin{defn} Let $(X,\pi)$ be a $K$-space. We define
\[
X_K = \{ T_\sigma   \;\mid \sigma\in K, \; T\in X, \; \sigma\leq \pi(T)\}
\]
\end{defn}

\begin{thm} Let $(X,\pi)$ be a $K$-space. 
Then $X_K$ is a ball complex.  Moreover
$X'$ is a simplicial subdivision of $X_K$. Also, $X_K$ is a subdivision of $X$. 

Let $f:(X,\pi_X)\to (Y,\pi_Y)$ is a map of $K$-spaces. The induced map $f':X'\to Y'$ of first derived 
complexes is then a map of ball complexes, $f_K: X_K\to Y_K$. 
\end{thm}
\begin{proof}  (The induced map $f'$ means the simplicial map $f':X'\to Y'$ for which $f'(bS)= b(f(S))$ for each $S\in X.$)
By Lemma \ref{ballem} the boundary of each $T_\sigma$ is a union of   balls of $X_K$ with smaller dimension  and
  \[
  T_\sigma ^\circ = \coprod \{A^\circ\;\mid A=\langle S_0, \dots, S_p\rangle\in D_\sigma T, \; A\notin \p^iD_\sigma T, \; A\notin \p^oD_\sigma T \} .
  \]
 This can be rewritten as:
  \begin{equation}\label{dst} T_\sigma^\circ = \coprod\{A^\circ\;\mid A=\langle S_0, \dots, S_p\rangle\in X', \; \sigma=\pi(S_p), \; T= S_0\},
\end{equation}

 By equation (\ref{dst}), for each $A\in X'$ there is a unique $T_\sigma\in X_K$ for which $A^\circ\subset T_\sigma ^\circ$.
 Therefore : $|X'|= \coprod \{T_\sigma ^\circ \; \mid \; \; T_\s\in X_K \}=|X_K|$.

 This proves that $X_K$ is a ball complex and that $X'$ is a subdivision of $X_K$. Because $T_\sigma \subset T$, we see $X_K$ is a subdivision of $X$.
 
 Now let $f:(X,\pi_X)\to (Y, \pi_Y) $ be a map of $K$-spaces. For each simplex $S\in X$  we see $f(S)\in Y$ because $f$ is simplicial. For each face $\s$ of 
 $\pi_X(S)$ in $K$, we see from the definitions  that  $f'(D_\s S)= D_\s f(S)$. So  $f'$ is a map of ball complexes, $f_K: X_K
 \to Y_K$.
\end{proof}

\section{The Isomorphism $\Phi_X:T\Delta^*X\cong C(X_K)$}\label{isom}
Our main theorem is:
\begin{thm} \label{MT}For each $K$-space $(X, \pi)$ the cellular chain complex
of $X_K$ with $R$ coefficients, denoted $C(X_K)$, comes with a natural $(R,K)$complex structure. There is defined (below) an isomorphism of $(R,K)$ chain complexes:
\[
\Phi_X :T\Delta^*X\cong C(X_K)\;.
\] 
 For each map $f:(X, \pi_X)\to (Y,\pi_Y)$ of $K$-spaces, the
square below commutes.
\[
\begin{CD}
T(\Delta^*X)@>T(f^*)>> T(\Delta^*Y)\\
@V\Phi_{X}VV                                @VV\Phi_{Y}V\\
C(X_K) @>f_K>>                       C(Y_K)
\end{CD} 
\]
\end{thm}

\

\begin{proof} 
Choose a basis $bK$ of oriented cells for $\Delta_*(K;R)$. Choose next, a basis $b_*X$ of oriented cells for  $\Delta_*(X;R)$. But choose the orientations in $b_*X$ so that
if $T\in b_*X$ and $\s\in bK$ are both $q$-cells, and if $\pi_*(T)=\pm \s\in \Delta_q(K;R)$, then:
\[
\pi_*(T) =(-1)^{dim\,\s}\s \in \Delta_q(K;R).
\]
We call such a pair, $(bK, b_*X)$ \emph{  an orientation for $(X,\pi)$.}

 Our first task is to construct the cellular chain complex $C_*(X_K;R)$ as the underlying $R$-complex of an $(R,K)$ complex $C(X_K)$ Define 
 \[
 C(X_K)= \Delta X\otimes_K\Delta^*K; \qquad C_*(X_K;R)=(\Delta X\otimes_K\Delta^*K)(K)
\]
For each oriented simplex $\rho$ of $K$ and oriented simplex $T$ of $X$, define
\[
[T_\rho] = T\otimes_K \rho^*\in C_{|T|-|\s|}(X_K;R)
\]
Define $bX_K= \{ [T_\rho]\mid  T\in b_*X,\; \rho \in bK, \;T_\rho\in X_K\}$. Then $bX_K$ is an
$R$-basis for $C_*(X_K;R)$ in bicorrespondence with the cells of $X_K$. Write $\p_q$
for the boundary map in $C_*(X_K;R)$, namely $(d^{\Delta X\otimes_K\Delta^*K})_q$\;.

But to justify these definitions, we must  check that $C_*(X_K;R)$ does  compute the cellular
homology or $X_K$. It  suffices to check, for any $[T_\rho]\in bX_K, $ that $\p_q([T_\rho])$ is a sum with $\pm1$  coefficients of those
$[S_\s]\in bX_K$ which are $(q-1)$-faces of $T_\rho$. (See \cite{cf}, for example.)

All proper faces of $T_\rho$ have the form $T_\s, $ for $\rho<\s$,\;  or $S_\rho$, for $S<T$.

Suppose $[T_\rho]\in bX_K$. So $T\in b_*X$, $\rho\in bK$. Set $\tau= \pi(T)\in K$.\; By (\ref{d_K}):
\begin{multline}\notag
\p_q[T_\rho]=  d^{\Delta X\otimes_K\Delta^*K} (T\otimes_K\rho^*)= 
\\\sum_{\{\s\mid \rho\leq \s\leq \tau\}} \{(d^{\Delta X}
(\s,\tau)T)\otimes\rho^*+(-1)^{|T|} T\otimes d^{\Delta^*K}(\s, \rho)\rho^*\}\\
=\sum_{S<T}[T,S][S_\rho] +(-1)^{1+|T_\rho|} \sum_{\rho<\s}[\s,\rho][T_\s]
\end{multline}
which is as required.

This completes the construction of the cellular chain complex of $X_K$, as
an $(R,K)$ complex, $C(X_K)$.

The $(R,K)$ isomorphism, $\Phi_X: T\Delta^* X\cong C(X_K)$ is simply:
\[
\Phi_X:= (\eps_{\Delta X}\otimes_K 1_{\Delta^*K}):T\Delta^*X=\Delta^{**}X \otimes_K\Delta^*K\lra \Delta X\otimes_K\Delta^*K= C(X_K).
\]
Naturality of $\Phi$ is obvious from the naturality of $\eps$.\end{proof}

 \section{ The McCrory Cap Product, \;$\Delta^*X$ and $\Delta X'$.}\label{sec 9}
 
 We now use the work of McCrory \cite{McC} to construct, for any $K$-space,
 $(X, \pi)$, an $(R,K)$ chain monomorphism serving two purposes:
 \[
 C(X_K)\overset{C_X}{\lra} \Delta X'.
 \]
  First, it defines an $(R,K)$ chain homotopy equivalence, $T\Delta^*X\simeq \Delta X'$.
  Second, $C_X$ identifies $C(X_K)$ with that $(R, K)$ subcomplex of $\Delta X'$
  which admits a basis consisting of one fundamental $q$-cycle, in
  $\Delta_q (D_\s T,\p D_\s T)\subset \Delta_q(X')$, for each $q$-cell $T_\s$ of $X_K$.
  (This will complete our geometric interpretation of $T$).

 \
 
 Let $K$ be a finite simplicial complex.   McCrory (see \cite{McC}, and also \cite{fl}) defines a map,
 $ c': \Delta_*(K;R)\otimes_R\Delta^*(K;R)\to \Delta_*(K';R)$ which he shows is chain homotopic to the composite,
 \[
  \Delta_*(K;R)\otimes_R\Delta^*(K;R)\overset{\cap}{\lra} \Delta_*(K;R)\overset{Sd}{\lra} \Delta_*(K';R) 
 \]
  where $\cap$ denotes the Whitney-Cech cap product. We will write $c_K$ for $c'$. We repeat his definition here with appropriate sign changes because McCrory's sign conventions differ slightly from ours. 
  
  For any $q$-simplex, $Q=\langle Q^0, Q^1, \dots Q^q\rangle \text{ of } K'$ \emph{in which  each $Q_i$ is oriented},  McCrory then defines
  \[ 
  \eps(Q) = [Q^0,Q^1][Q^1,Q^2]\dots [Q^{q-1},Q^q]. 
  \] 
  This is independent of the orientations on $Q_1, Q_2,\dots Q_{q-1}$. If $q=0$, set $\eps(Q)=0$.
  
  For any $n$- simplex $\tau$ and $(n-q)$-simplex $\s$ of $K$, each simplex
 $Q=\langle Q^0, Q^1, \dots Q^q\rangle$ of $D(\s, \tau)_q$ satisfies: $Q_0= \tau; \;Q_q=\s$. Therefore, $\eps(Q)$ makes sense if $\tau$ and $\s$ are oriented simplices 
 chosen from some basis $bK$ of oriented simplices for $\Delta_*(K;R)$ (but not if $\s=-\tau).$
  
  The \emph{McCrory Cap Product},   $\Delta_*(K;R)\otimes_R\Delta^*(K;R)\overset{c_K}{\lra} \Delta_*(K';R)$  is the map defined by:
  \[
  \boxed{
  c_K(\tau\otimes \s^*)= \sum_{Q\in D(\s,\tau)_q}(-1)^{dim (\s)} \eps(Q)Q
  }
  \]
  for any oriented simplices $\s, \tau$ in some basis $bK$.  Here $q= dim(\tau)-dim(\s)$. Note this is zero unless $\s\leq \tau$. Note that $c_K$ does not change
  if we change the basis.
  
  $c_K$ is a chain map. We reprove this in the Appendix, Section \ref{App}, because of the sign changes and because McCrory's proof, \cite{McC} p.155 lines 7-8, is only a sketch.
  
  \
  
  Now suppose $(X,\pi)$ is a $K$-space.  
  
  Note that if $T$ and $\s$ are oriented simplices of $X$ and  $K$ and $q=dim(T)-dim(\s)\neq 0$:
  \[
  c_X(T\otimes_R \pi^*\s^*)=\sum_{Q\in (D_\s, T)_q}(-1)^{dim (\s)} \eps(Q)Q\;\;\in \Delta_qX'(\s)
  \]
  (because $D_\s T=\cup\{D(S,T)\; |\;  S\in X, dim(S)=dim(\s), \pi(S)= \s \}$).
  This formula still makes sense and is true  if $q=0$ and $\pi_*(T)\neq -\s$).
  
  In this way, $c_X\circ(1\otimes \pi^*)$ defines an $(R,K)$ chain map,
  \[
  c_X\circ(1\otimes \pi^*): \Delta X\otimes_R\Delta^*K\lra \Delta X'
  \]
  
  \begin{prop} There is a unique $(R,K)$ chain map 
  \[
    C_X: C(X_K)=\Delta X\otimes_K \Delta^*K \lra \Delta X'
    \]
    satisfying:
    \[
    c_X\circ(1\otimes \pi^*)=C_X\circ \pi_{\Delta X, \Delta^*K} 
    \]
    $C_X$ is an $(R,K)$  monomorphism. For all $q$-cells $T_\s$ of
    $X_K$, with $q\neq 0$,
  \[
  C_X(T\otimes_K \s^*)= \sum_{Q\in (D_\s T)_q}(-1)^{dim (\s)} \eps(Q)Q,  
  \]
   For a $0$-cell $T_\s$, of $X_K$, with $T\in\Delta_n(X;R), \s\in\Delta_n(K;R)$ are
  oriented so that $\pi_*(T)=\s$, then
  \[
  C_X(T\otimes_K\s^*)= (-1)^{dim; T}\langle bT\rangle, \quad (bT \text{is the barycenter of } T). 
  \]
  \end{prop}
  \begin{proof}
  Note that $c_X(T\otimes_R \pi^*\s^*)=0$ unless $\pi(T)\geq \s$.  Also $c_X(T\otimes_R \pi^*\s^*)\in \Delta X'(\s)$ for all $\s\in K$ and $T\in X$ because each $q$-cell $Q\in D_\s T$ is in $\Delta_qX'(\s)$ if $q= dim(T_\s)$.
  
  So $c_X\,\circ(1\otimes \pi^*):\Delta X\otimes_R \Delta^*K\to \Delta X'$  is an $(R,K)$ chain map annihilating each 
  $\Delta X(K-st(\s))\otimes_R \Delta ^*K(\s)$.  Hence  there is a unique $(R,K)$ chain map monomorphism,
  $ \Delta X\otimes_K \Delta^*K\overset{C_X}{\lra} \Delta X'$ such that
  $c_X\circ(1\otimes \pi^*) = C_X\circ \pi_{\Delta X, \Delta^*K} $ .   The calculation follows if $q\neq 0$. If $q=0$, then $(\pi_{\mid T})^* \s^*=T^*$, so
  
  \[
  C_X(T\otimes_K \s^*)= c_X(T\otimes T^*) =(-1)^{dim\,T}\sum_{Q\in D(T,T)_0}Q= (-1)^{dim\;T}\langle bT\rangle
  \]
  
  Clearly $C_X$ is natural in $(X,\pi)$.
    \end{proof}
  \begin{rem}: If we choose an orientation $(bK, b_*X)$ for $X_K$, then for each $0$-cell
  $T_\s = T_{\pi(T)} $ of $X_K$ with $[T_\s]\in bX_K$, we have $C_X([T_\s])= \langle bT\rangle$.
  \end{rem}

    \begin{cor} \label{fnc}  For each $q$-cell $T_\s$ of $X_K$, $C_X(T\otimes \s^*)$
  is a fundamental cycle, in $\Delta_q(D_\s T, \p D_\s T)$ for the $q$-manifold $D_\s T$.
     \end{cor}
     \begin{proof}$C_X(T\otimes_K \s^*)$ is  a fundamental cyle in $\Delta_q(D_\s T, \p D_\s T)$ since $C_X$ is a chain map and since each 
    $ Q\in (D_\s T)_q $ appears with coefficient $\pm 1$ in $C_X(T\otimes_K \s^*)$.
     \end{proof}
     \begin{thm} \label{che} For each $K$-space $(X,\pi)$, the map $C(X_K)\overset{C_X}{\lra} \Delta X'$ is an $(R,K)$ chain homotopy equivalence.
     \end{thm} 
     \begin{proof} By \ref{fnc}, for all $T_\s$ , $C_X$ restricts to a homotopy equivalence, 
     \[
     C_*(T_\s,\p T_\s; R)\to \Delta_*(D_\s(T),\p D_\s(T);R)
     \]
       and it takes chains on any subcomplex of $X_K$ to chains on its subdivision. By an induction-excision argument on the number of cells in the subcomplex one sees $C_X$ yields a homology equivalence and then a chain homotopy equivalence on each such subcomplex. So $C_X(\s,\s) $ is an $R$-chain equivalence for each $\s$. Therefore $C_X $ is an $(R,K)$ chain equivalence.
     \end{proof}  
     Together,  \ref{che} and \ref{MT} clearly prove:
     \begin{cor}
    $T\Delta^*X\overset{ C_X \Phi_X}{\lra} \Delta X'$ is an $(R,K)$ chain homotopy equivalence.\quad
    Consequently $e_{\Delta^*X}\circ T(C_X\Phi_X)$ is an explicit $(R,K)$ chain
    homotopy equivalence,  
    \[
    T\Delta X'\simeq \Delta^*X.
    \]
     \end{cor} 
   
 \
 
\section {Appendix}\label{App}
We must prove:
\begin{prop}$\Delta_*(K;R)\otimes_R\Delta^*(K;R)\overset{c_K}{\lra} \Delta_*(K';R)$ is a chain map. \newline That is
to say, for any oriented simplices $\s, \tau$ in some basis $bK$   for $\Delta K$, with \newline$ p= dim(\tau)-dim(\s)$, 
\[
d^{K'} c_K(\tau\otimes \s^*) = c_K\{d^K \tau \otimes \s^*+(-1)^{dim(\tau)} \tau\otimes d^{\Delta^*(K)}\s^*\} 
\]
\end{prop}
where, by the definitions, 
\[d^K\tau=\sum_{\rho\in b_K}[\tau,\rho] \rho,\qquad \quad d^{\Delta^*(K)}\s^*= (-1)^{dim(\s)+1}\sum_{\rho\in b_K}[\rho, \s] \rho^*
\]
and for any $p$-simplex $Q=\langle Q^0, Q^1,\dots Q^p\rangle$ of $K'$, 
\[
d^{K'} Q= \sum_{i=0}^p (-1)^i d^i(Q); \qquad  d^i(Q)= \langle Q^0, Q^1,\dots \hat{Q^i}\dots Q^p\rangle
\]
\begin{proof} We first prove:\quad $d^oc(\tau\otimes\s^*) = c(d^K\tau \otimes \s^*)$, where $c=c_K$.
\begin{multline}\notag
d^0c(\tau\otimes \s)= (-1)^{dim\s}\sum_{Q\in D(\s,\tau)_p}\eps(Q)\langle Q^1,\dots Q^p\rangle = \\(-1)^{dim\s}\sum_{\rho\in b_K} [\tau,\rho]\sum_{P\in D(\s, \rho)} \eps(P)P =c(\sum_{\rho\in b_K} [\tau,\rho] \rho\otimes \s^*) = c(d^K\tau\otimes \s^*).
\end{multline}
Next we show: $(-1)^pd^p c(\tau\otimes\s^*)=(-1)^{dim(\tau)}c( \tau\otimes d^{\Delta^*(K)}\s^*) $:

\begin{multline}
(-1)^pd^p c(\tau\otimes\s^*)=(-1)^{p+dim\s}\sum_{Q\in D(\s,\tau)_p}\eps(Q)\langle \tau, Q^1\dots Q^{p-1}\rangle=\\
(-1)^{p+1} c(\tau\otimes \sum_{\rho\in b_K} [\rho,\s]\rho^*)=(-1)^{dim(\tau)}c(\tau\otimes d^{\Delta^*(K)}\s^*)
\end{multline}
Finally we prove $d^ic(\tau\otimes \s^*)=0$ for $0<i<p$.

For such $i$ and for $Q\in D(\s, \tau)$ note $d^iQ=\langle \tau,\dots \s\rangle\in D(\s,\tau)- \p D(\s,\tau)$.
So suppose $P$ is a $p-1$ simplex of the form $d^iQ$ in the $p$ manifold $D(\s,\tau)$. Then there is exactly one other $S\in D(\s,\tau)_p$ having $Q$ as a face.
We can identify $S$ by listing the vertices of $\tau$ as $v_0, \dots v_n$ so that
$Q^j =[ v_j, \dots v_n]$ for all $j$. Define $S^i=[v_0\dots v_{i-1}, v_{i+1}\dots v_n]$ and define
$S^j=Q^j$ for $j\neq i$. Then $S:= \langle S^0, S^1, \dots S^p\rangle$ in $ D(\s,\tau)_p$ satisfies $d^iS=P;\; \eps(S)=-\eps(Q)$ so $P$ must appear with zero coefficient in $d^i c(\tau\otimes \s^*)$ for all $p-1$ simplices $P$. So $ d^ic(\tau\otimes \s^*)=0$.\end{proof}

\end{document}